\documentclass[11pt]{amsart}
\usepackage{amsfonts}
\usepackage{mathrsfs}
\usepackage{amssymb,latexsym}
\usepackage{pstcol,pstricks,color}

 \numberwithin{equation}{section}

\newtheorem{theorem}{Theorem}[section]

\def\qed{\hfill $\Box$}
\def\pf{\noindent {\it Proof.} }

\title{A simple bijection between binary trees and colored ternary trees}

\begin{document}
\maketitle

\begin{center}
Yidong Sun

Department of Mathematics, Dalian Maritime University, 116026 Dalian, P.R. China\\[5pt]

{\it sydmath@yahoo.com.cn}
\end{center}\vskip0.5cm

\subsection*{Abstract} In this short note, we first present a simple
bijection between binary trees and colored ternary trees and then
derive a new identity related to generalized Catalan numbers.

\medskip

{\bf Keywords}: Binary tree; Ternary tree; Generalized Catalan
number

\noindent {\sc 2000 Mathematics Subject Classification}: Primary
05C05; Secondary05A19

\section{Introduction}

Recently, Mansour and the author \cite{mansun} obtained an identity
involving $2$-Catalan $C_{n,2}=\frac{1}{2n+1}\binom{2n+1}{n}$
numbers and $3$-Catalan numbers
$C_{n,3}=\frac{1}{3n+1}\binom{3n+1}{n}$, i.e.,
\begin{eqnarray} \label{eqn 1}
\sum_{p=0}^{[n/2]}\frac{1}{3p+1}\binom{3p+1}{p}\binom{n+p}{3p}&=&\frac{1}{n+1}\binom{2n}{n}.
\end{eqnarray}
In this short note, we first present a simple bijection between
complete binary trees and colored complete ternary trees and then
derive a general identity, i.e.,
\begin{theorem}\label{theo 1}
For any integers $n,m\geq 0$, there holds
\begin{eqnarray}\label{eqn 2}
\sum_{p=0}^{[n/2]}\frac{m}{3p+m}\binom{3p+m}{p}\binom{n+p+m-1}{n-2p}
&=&\frac{m}{2n+m}\binom{2n+m}{n}.
\end{eqnarray}
\end{theorem}
\vskip0.3cm

\section{A bijective algorithm for binary and ternary trees}

A {\em colored ternary trees} is a complete ternary tree such that
all its vertices are signed a nonnegative integer called {\em color
number}. Let $\mathscr{T}_{n,p}$ denote the set of colored ternary
trees $T$ with $p$ internal vertices such that the sum of all the
color numbers of $T$ is $n-2p$. Define $\mathscr{T}_n=\bigcup_{p=
0}^{[n/2]}\mathscr{T}_{n,p}$. Let $\mathscr{B}_{n}$ denote the set
of complete binary trees with $n$ internal vertices. For any $B\in
\mathscr{B}_n$, let $P=v_1v_2\cdots v_k$ be a path of {\em length}
$k$ of $B$ (viewing from the root of $B$). $P$ is called a {\em
R-path}, if (1) $v_i$ is the right child of $v_{i-1}$ for $2\leq
i\leq k$ and (2) the left child of $v_i$ is a leaf for $1\leq i\leq
k$. In addition, $P$ is called a {\em maximal R-path} if there
exists no vertex $u$ such that $uP$ or $Pu$ forms a $R$-path. $P$ is
called an {\em L-path}, if $v_i$ is the left child of $v_{i-1}$ for
$2\leq i\leq k$. $P$ is called a {\em maximal L-path} if there
exists no vertex $u$ such that $uP$ or $Pu$ forms an $L$-path.

Note that the definition of $L$-path is different from that of
$R$-path. Hence, if $P$ is a maximal $R$-path, then (1) the right
child $u$ of $v_k$ must be a leaf or the left child of $u$ is not a
leaf; (2) $v_1$ must be a left child of its father (if exists) or
the father of $v_1$ has a left child which is not a leaf. If $P$ is
a maximal $L$-path, then (1) $v_k$ must be a leaf which is also a
left child of $v_{k-1}$; (2) $v_1$ must be the right child of its
father (if exists).

\begin{theorem}\label{theo 2}
There exists a simple bijection $\phi$ between $\mathscr{B}_{n}$ and
$\mathscr{T}_{n}$.
\end{theorem}
\pf We first give the procedure to construct a complete binary tree
from a colored complete ternary tree.
\begin{itemize}
\item[Step 1.]  For each vertex $v$ of $T\in \mathscr{T}_{n}$ with
color number $c_v=k$, remove the color number and add a $R$-path
$P=v_1v_2\cdots v_k$ of length $k$ to $v$ such that $v$ is a right
child of $v_k$ and $v_1$ is a child of the father (if exists) of
$v$, and then annex a left leaf to $v_i$ for $1\leq i\leq k$. See
Figure \ref{fDD1}(a) for example.
\begin{figure}[h]
\setlength{\unitlength}{0.4mm}
\begin{center}
\begin{pspicture}(11,5.2)
\psset{xunit=19pt,yunit=19pt}\psgrid[subgriddiv=1,griddots=5,
gridlabels=4pt](0,0)(17,8)

\psline(2,6.5)(2,5)(.5,3)\psline(2,5)(2,3)\psline(2,5)(3.5,3)

\pscircle*(2,5){0.075}\pscircle*(2,3){0.075}
\pscircle*(.5,3){0.075}\pscircle*(3.5,3){0.075}

\put(.9,3.4){$v$}\put(1.7,3.5){$\tiny{c_v=2}$}\put(.2,1.6){$\tiny
T_1$}\put(1.2,1.6){$\tiny T_2$}\put(2.2,1.6){$\tiny{T_3}$}
\put(3,2.9){$\Longleftrightarrow$}\put(3.1,.15){$(a)$}

\psline(7,7.8)(7,6.5)(6,5)\psline(7,6.5)(8,5)(7,3.5)\psline(8,5)(9,3.5)(8,2)
\psline(9,3.5)(9,2)\psline(9,3.5)(10,2)

\pscircle*(7,6.5){0.075}\pscircle*(6,5){0.075}
\pscircle*(8,5){0.075}\pscircle*(7,3.5){0.075}\pscircle*(9,3.5){0.075}
\pscircle*(8,2){0.075}\pscircle*(9,2){0.075}\pscircle*(10,2){0.075}

\put(6.2,2.3){$v$}\put(4.8,4.3){$v_1$}\put(5.5,3.3){$v_2$}\put(5.2,.9){$\tiny
T_1$}\put(5.9,.9){$\tiny T_2$}\put(6.6,.9){$\tiny{T_3}$}
\put(7,2.9){$\Longleftrightarrow$}\put(7.1,.15){$(b)$}

\psline(13,7.8)(13,6.5)(12,5)\psline(13,6.5)(14,5)(13,3.5)\psline(14,5)(15,3.5)(14,2)(13,.5)
\psline(15,3.5)(16,2)\psline(14,2)(15,.5)

\pscircle*(13,6.5){0.075}\pscircle*(12,5){0.075}
\pscircle*(14,5){0.075}\pscircle*(13,3.5){0.075}\pscircle*(15,3.5){0.075}
\pscircle*(14,2){0.075}\pscircle*(13,.5){0.075}\pscircle*(15,.5){0.075}
\pscircle*(16,2){0.075}

\put(10.2,2.3){$v$}\put(9.6,1.3){$v'$}\put(8.8,4.3){$v_1$}\put(9.5,3.3){$v_2$}\put(8.9,.2){$\tiny
T_1$}\put(10.2,.2){$\tiny T_2$}\put(10.6,.9){$\tiny{T_3}$}

\end{pspicture}
\caption{\ }\label{fDD1}
\end{center}
\end{figure}

{\item[Step 2.] Let $T^*$ be the tree obtained from $T$ by Step $1$.
For any internal vertex $v$ of $T^*$ which has out-degree $3$, let
$T_1,T_2$ and $T_3$ be the three subtree of $v$. Remove the subtree
$T_1$ and $T_2$ , annex a left child $v'$ to $v$ and take $T_1$ and
$T_2$ as the left and right subtree of $v'$ respectively. See Figure
\ref{fDD1}(b) for example.}

\end{itemize}

It is clear that any $T\in \mathscr{T}_n$, after Step $1$ and $2$,
generates a binary tree $B\in \mathscr{B}_n$.

Conversely, we can obtain a colored ternary tree from a complete
binary tree as follows.
\begin{itemize}
\item[Step 3.] Choose any
maximal $L$-path of $B\in \mathscr{B}_n$ of length $k\geq 3$, say
$P=v_1v_2\cdots v_k$, then each $v_{2i-1}$ absorbs its left child
$v_{2i}$ for $1\leq i\leq [(k-1)/2]$. See Figure \ref{fDD2}(a) for
example.
\begin{figure}[h]
\setlength{\unitlength}{0.4mm}
\begin{center}
\begin{pspicture}(11,5.2)
\psset{xunit=19pt,yunit=19pt}\psgrid[subgriddiv=1,griddots=5,
gridlabels=4pt](-2,0)(18,8)

\psline(-.5,6.5)(.5,8)(1.5,6.5)\psline[linewidth=2pt](1.5,6.5)(.5,5)(-1.5,2)\psline(1.5,6.5)(2.5,5)(1.7,3.5)
\psline[linewidth=2pt](3.5,3.5)(2.5,2)(1.7,.5)
\psline(.5,5)(1.2,3.5)\psline(-.5,3.5)(.5,2)(1.2,.5)
\psline(.5,2)(-.5,.5)\psline(2.5,5)(3.5,3.5)(4.5,2)\psline(2.5,2)(3.5,.5)

\pscircle*(-.5,6.5){0.075}\pscircle*(.5,8){0.075}
\pscircle*(1.5,6.5){0.075}\pscircle*(.5,5){0.075}\pscircle*(-.5,.5){0.075}\pscircle*(1.2,.5){0.075}
\pscircle*(2.5,5){0.075}\pscircle*(1.7,3.5){0.075}\pscircle*(3.5,3.5){0.075}
\pscircle*(2.5,2){0.075}\pscircle*(1.7,.5){0.075}\pscircle*(3.5,.5){0.075}
\pscircle*(4.5,2){0.075}\pscircle*(-1.5,2){0.075}\pscircle*(.5,2){0.075}
\pscircle*(-.5,3.5){0.075}\pscircle*(1.2,3.5){0.075}

\put(3,2.9){$\Longleftrightarrow$}\put(3.1,.15){$(a)$}

\psline(6,6.5)(7,8)(8,6.5)(7,5)(6,3.5)\psline[linewidth=2pt](7,5)(7.5,3.5)
\psline(7.5,3.5)(8,2)
\psline(7.5,3.5)(7,2)\psline(8,6.5)(8,5)
\psline(8,6.5)(9,5)(8.5,3.5)\psline(9,5)(10,3.5)(9,2)
\psline(10,3.5)(10,2)\psline(10,3.5)(11,2)

\pscircle*(8,6.5){0.075}\pscircle*(7,5){0.1}\pscircle*(7,8){0.1}
\pscircle*(6,6.5){0.075}\pscircle*(8,5){0.075}\pscircle*(7,2){0.075}\pscircle*(8,2){0.075}
\pscircle*(9,5){0.1}\pscircle*(8.5,3.5){0.075}\pscircle*(7.5,3.5){0.1}
\pscircle*(6,3.5){0.075}\pscircle*(10,3.5){0.075}
\pscircle*(9,2){0.075}\pscircle*(10,2){0.075}\pscircle*(11,2){0.075}

\put(7.7,2.9){$\Longleftrightarrow$}\put(7.8,.15){$(b)$}

\psline(15,6.5)(14,5)\psline(15,6.5)(15,5)\psline(15,6.5)(16,5)(15,3.5)
\psline(16,5)(17,3.5) \psline(16,5)(16,3.5)

\pscircle*(15,6.5){0.075}\pscircle*(15,5){0.075}\pscircle*(14,5){0.075}
\pscircle*(16,5){0.075}\pscircle*(15,3.5){0.075}\pscircle*(17,3.5){0.075}
\pscircle*(16,3.5){0.075}

\put(10.2,4.3){$1$}\put(10.9,3.3){$1$}\put(9.2,2.9){$2$}\put(9.9,2.9){$0$}
\put(9.9,1.9){$0$}\put(10.6,1.9){$0$}\put(11.4,1.9){$0$}

\end{pspicture}
\caption{\ }\label{fDD2}
\end{center}
\end{figure}

\item[Step 4.] Choose any maximal $R$-path of $T'$ derived from $B$
by Step $3$, say $Q=u_1u_2\cdots u_k$, let $u$ be the right child of
$u_k$, then $u$ absorbs all the vertex $u_{1},u_2,\dots,u_k$ and
assign the color number $c_u=k$ to $u$. See Figure \ref{fDD2}(b) for
example. Hence we get a colored ternary tree. \qed

\end{itemize}

Given a complete ternary tree $T$ with $p$ internal vertices, there
are totally $3p+1$ vertices, choose $n-2p$ vertices repeatedly,
define the color number of a vertex to be the times of being chosen.
Then there are $\binom{n+p}{n-2p}$ colored ternary trees in
$\mathscr{T}_{n}$ generated by $T$. Note that
$\frac{1}{3p+1}\binom{3p+1}{p}$ and $\frac{1}{2n+1}\binom{2n+1}{n}$
count the number of complete ternary trees with $p$ internal
vertices and complete binary trees with $n$ internal vertices
respectively \cite{stanley}. Then the bijection $\phi$ immediately
leads to (\ref{eqn 1}).

To prove (\ref{eqn 2}), consider the forest of colored ternary trees
$F=(T_1,T_2,\dots, T_m)$ with $T_i\in \mathscr{T}_{n_i}$ and
$n_1+n_2+\cdots+n_m=n$, define $\phi(F)=(\phi(T_1),\phi(T_2),\dots,
\phi(T_m))$, then it is clear that $\phi$ is a bijection between
forests of colored ternary trees and forests of complete binary
trees. Note that there are totally $m+3p$ vertices in a forest $F$
of complete ternary trees with $m$ components and $p$ internal
vertices, so there are $\binom{m+n+p-1}{n-2p}$ forests of colored
ternary trees with $m$ components, $p$ internal vertices and the sum
of color numbers equal to $n-2p$. It is clear  \cite{stanley} that
$\frac{m}{3p+m}\binom{3p+m}{p}$ counts the number of forests of
complete ternary trees with $p$ internal vertices and $m$
components, and that $\frac{m}{2n+m}\binom{2n+m}{n}$ counts the
number forests of complete binary trees with $n$ internal vertices
and $m$ components. Then the above bijection $\phi$ immediately
leads to (\ref{eqn 2}).

\section{Further comments}

It is well known \cite{stanley} that the $k$-Catalan number
$C_{n,k}=\frac{1}{kn+1}\binom{kn+1}{n}$ counts the number of
complete $k$-ary trees with $n$ internal vertices, whose generating
function $C_k(x)$ satisfies
\begin{eqnarray*}
C_k(x)=1+xC_k^{k}(x).
\end{eqnarray*}

Let $G(x)=\frac{1}{1-x}C_3(\frac{x^2}{(1-x)^3})$, then one can
deduce that
\begin{eqnarray*}
G(x)&=&\frac{1}{1-x}C_3(\frac{x^2}{(1-x)^3})\\
    &=&\frac{1}{1-x}(1+\frac{x^2}{(1-x)^3}C_3^3(\frac{x^2}{(1-x)^3}))\\
    &=&\frac{1}{1-x}(1+x^2G^3(x)),
\end{eqnarray*}
which generates that $G(x)=C_2(x)$ which is the generating function
for Catalan numbers.

By Lagrange inversion formula, we have
\begin{eqnarray*}
C_3^m(x)&=&\sum_{p\geq 0}\frac{m}{3p+m}\binom{3p+m}{p}x^p, \\
C_2^m(x)&=&\sum_{n\geq 0}\frac{m}{2n+m}\binom{2n+m}{n}x^n.
\end{eqnarray*}
Then
\begin{eqnarray*}
G^m(x)&=&\sum_{p\geq
0}\frac{m}{3p+m}\binom{3p+m}{p}\frac{x^{2p}}{(1-x)^{3p+m}} \\
&=&\sum_{n\geq
0}x^n\sum_{p=0}^{[n/2]}\frac{m}{3p+m}\binom{3p+m}{p}\binom{n+p+m-1}{n-2p}.
\end{eqnarray*}
Comparing the coefficient of $x^n$ in $C_2^m(x)$ and $G^m(x)$, one
obtains Theorem \ref{theo 2}.

Similarly, let $F(x)=\frac{1}{1-x}C_k(\frac{x^{k-1}}{(1-x)^{k}})$,
then $F(x)=\frac{1+xF(x)}{1-x^{k-1}F^{k-1}(x)}$, using Lagrange
inversion formula for the case $k=5$, one has
\begin{eqnarray}\label{eqn 3.1}
\lefteqn{\sum_{p=0}^{[n/4]}\frac{m}{5p+m}\binom{5p+m}{p}\binom{n+p+m-1}{n-4p}}\\
\nonumber
&=&\sum_{p=0}^{[n/2]}(-1)^p\frac{m}{m+n}\binom{m+n+p-1}{p}\binom{m+2n-2p-1}{n-2p},
\end{eqnarray}
which, in the case $m=1$, leads to
\begin{eqnarray}\label{eqn 3.2}
\sum_{p=0}^{[n/4]}\frac{1}{4p+1}\binom{5p}{p}\binom{n+p}{5p}
&=&\sum_{p=0}^{[n/2]}(-1)^p\frac{1}{n+1}\binom{n+p}{n}\binom{2n-2p}{n}.
\end{eqnarray}

One can be asked to give a combinatorial proof of (\ref{eqn 3.1}) or
(\ref{eqn 3.2}).

\vskip1cm

\section*{Acknowledgements} The authors are grateful to the
anonymous referees for the helpful suggestions and comments.


\end{document}